\documentclass[10pt, color, amstex]{amsart}
\newcommand{\pof}{\noindent{\em Proof: }}
\usepackage{color}
\newcommand{\s}[1]{\mathcal{#1}}

\newcommand{\tr}{\operatorname{tr}}

\newcommand{\Deg}{\operatorname{deg}}
\newcommand{\str}{\operatorname{str}}

\newtheorem{Thm}{Theorem}[section]
\newtheorem{Def}[Thm]{Definition} \newtheorem{Rem}[Thm]{Remark}
\newtheorem{Lem}[Thm]{Lemma} 
\newtheorem{Con}[Thm]{Conjecture}
\newtheorem{Prop}[Thm]{Proposition}
\numberwithin{equation}{section}

\begin{document}

\fontsize{.5cm}{.5cm}\selectfont\sf

\vskip0.5cm
\medskip

\title{The exponential nature and positivity}
\author{Hans Plesner jakobsen, Hechun Zhang }\thanks{This research was
  supported in part by Project AM14:21/ NPP33 under the Sino-Danish
  Scientific and Technological Cooperation. The second author was also
  supported by the Chinese National Science Foundation; project number 10471070}

\address{
  Institute of Mathematics\\ University of Copenhagen\\Universitetsparken 5\\
   DK-2100, Copenhagen,
  Denmark\\Department of Mathematical Sciences, Tsinghua University,
  Beijing, 100084, P. R. China} \email{jakobsen@@math.ku.dk,
hzhang@@math.tsinghua.edu.cn}\date{\today}

\begin{abstract}
  In the present article, a basis of the coordinate algebra of the
  multi-parameter quantized matrix is constructed by using
    an elementary method due to Lusztig. The construction depends
  heavily on an anti-automorphism, the bar action. The exponential
  nature of the bar action is derived which provides an inductive way
  to compute the basis elements. By embedding the basis into the dual
  basis of Lusztig's canonical basis of $U_q(n^-)$, the positivity
  properties of the basis as well as the positivity properties of the
  canonical basis of the modified quantum enveloping algebra of type
  $A$, which has been conjectured by Lusztig, are proved.
\end{abstract}

 \maketitle

\section{introduction}

The coordinate algebra of the multi-parameter quantized matrix has
been introduced by Artin, Schelter and Tate \cite{ast}. In \cite{jz2},
Jakobsen and Zhang introduced a class of quadratic algebras which are
not bi-algebras in general, but are quite similar to the coordinate
algebra of the quantum matrix with one-parameter in the representation
aspect. In the present paper, we construct a class of quadratic
algebras ${\s O}_{q, P, Q}(M(n))$ which includes the above two classes
of algebras as special cases, by using bi-character deformation. The
algebra ${\s O}_{q, P, Q}(M(n))$ is not a bi-algebra in general.
However, since we only consider the basis of the coordinate algebra,
we need not to deal with the multiplication of the matrices. We will
show that our algebra is a bi-character deformation of the quantum
matrix algebra algebra considered in \cite{pw}. Hence, the quantum
determinant and quantum minors can be defined in the exactly the same
way as in \cite{pw} and all of the properties of the quantum minors of
the so-called official quantum matrix algebras transfer to the minors
of the present algebra after slight modifications.  By using a method
by Lusztig, we construct a nice basis of the algebra -- the dual
canonical basis -- which is "invariant`` under the multiplication of
certain quantum minors. The main ingredient of the construction is a
nice anti-automorphism, the bar action. We prove that the bar action
has an exponential nature by introducing some simple operators
$T_{ij}^{st}, \overline{T_{ij}^{st}}$. Hence, our construction
provides an inductive algorithm for computing the basis elements. We
then compute the bases for certain special cases, and propose a
conjecture for the general case. Embedding the algebra $O_q(M(n))$
into the negative part $U_q(A_{2n-1})^-$ of the quantum enveloping
algebra $U_q(A_{2n-1})$, we show that our dual canonical basis is a
subset of the dual canonical basis (after a slight modification) of
Lusztig's canonical basis; this provides an interpretation of the
coefficients of the expansion of our dual canonical basis elements in
terms of the modified monomials $Z(A)$ and enables us to prove the
positivity properties of our basis and, by duality, the positivity
property of the canonical basis of the modified quantum enveloping
algebra of type $A$ is proved which was conjectured by Lusztig
\cite{lu2}.

 \medskip

\section{Lusztig's construction of the basis}

\medskip

Our method to construct the basis is a modification of Lusztig's
construction in \cite{lu1}, see also \cite{du}. We refer
  to this construction as Lusztig's elementary method.

Let $\Gamma$ be an abelian group with a total ordering which is
compatible with the group structure on $\Gamma$. Let $\Gamma_+$ be
the set of elements of $\Gamma$ which are strictly positive for
this ordering and let $\Gamma_-=(\Gamma_+)^{-1}$. Let $a\mapsto
\bar{a}$ be the involution of the group ring ${\mathbb Z}[\Gamma]$
which takes $\gamma$ to $\gamma^{-1}$. Let $V$ be a free ${\mathbb
Z}[\Gamma]$ module with a basis $\{t_i|i\in I\}$, where the index
set $I$ has an ordering $\le$.  Assume that there is a map
$$-: V\longrightarrow V,$$
satisfying $ \overline{a.v}=\bar{a}\bar{v}$ for all $a\in {\mathbb
Z}[\Gamma]$ and $v\in V$. Furthermore, assume that
$$\bar{t_i}=\sum a_{ij}t_j$$
with $a_{ii}=1$ and $a_{ij}\ne 0$ only if $j\le i$. In \cite{lu1},
Lusztig proved that
\begin{Prop} Given $i\in I$, there is a unique element $b_i\in V$
such that
$$\overline{b_i}=b_i,$$
and
$$b_i=\sum_{j\le i} h_{ij}b_j,$$
where $h_{ii}=1$ and for any $j<i$, $h_{ij}\in {\mathbb
Z}[\Gamma^+]$. The elements $b_i$ form a basis of $V$.\end{Prop}

The coefficients $h_{ij}$ satisfy a system of equations
\begin{eqnarray}h_{ii}&=&1,\\\nonumber
\overline{h_{ij}}-h_{ij}&=&\sum_{i<k<j}a_{ik}h_{kj}.\end{eqnarray}
Hence, the coefficients $h_{ij}$ can be computed inductively,
provided the $a_{ij}$ are known.

\section{  The bi-character deformation and the multi-parameter  quantum  matrix space}

Let $P=(p_{ij})$ and $Q=(q_{ij})$  be matrices whose entries
satisfy
$$p_{ii}=q_{ii}=p_{ij}p_{ji}=q_{ij}q_{ji}=1,$$
where $p_{ij},q_{ij}$ ($1\le i<j\le n)$, and, later, $q$, are
independent
 variables. The base field in the rest of the paper is
  $K={\mathbb Q}(q, p_{ij},q_{ij}|i<j)$.

  Let $G$ be an abelian semi-group. A semi-group homomorphism from
  $G\times G$ to  the multiplicative group $K^*$ is
  called a bi-character of $G$.

Let $A$ be an associative algebra with a $G\times G$ gradation:
$$A=\oplus_{g,h\in G} A_{g,h},$$
satisfying
$$ A_{g_1, h_1}A_{g_2, h_2}\subset A_{g_1+g_2, h_1+h_2}$$
for all $g_1,g_2, h_1, h_2\in G$.

Let $\phi$ and $\psi$ be two bi-characters of $G$. For $a\in
A_{g_1, h_1}$ and $b\in A_{g_2, h_2}$,  one may define a new
multiplication

$$a*b=\phi(g_1,g_2)\psi(h_1,h_2)ab.$$

It is easy to check that the new multiplication $*$ is associative
and this defines a new  associative algebra  structure on $A$
which  is called a bi-character deformation of the algebra $A$.

Define the coordinate algebra ${\s O}_{q,P,Q}(M(n))$ to be
  the
associative algebra generated by $n^2$ generators $Z_{ij}$ subject
to the defining relations:
\begin{eqnarray}Z_{st}Z_{ij}&=& p_{si}^2q_{tj}^2Z_{ij}Z_{st}+(q^2-1)p_{si}^2Z_{it}Z_{sj},
\text{ if }s>i,t>j,\\\nonumber
Z_{st}Z_{ij}&=&q^2p_{si}^2q_{jt}^{-2}Z_{ij}Z_{st}, \text {if
}s>i,t\le j,\\\nonumber
Z_{it}Z_{ij}&=&q_{tj}^2Z_{ij}Z_{it},\end{eqnarray}

\begin{Rem}
We get the official $2\times2$ matrix algebra with the following
choices: $p_{21}=q^{-1/2}\;;\ q_{21}=q^{1/2}$ (with ``$q^{-1}$ relations'').
\end{Rem}

\begin{Prop} The algebra ${\s O}_{q,P,Q}(M(n))$ is a bi-character
deformation of the coordinate algebra of the quantum matrix space
of Dipper-Donkin (\cite{dd}).\end{Prop}

\pof The coordinate algebra of the quantum matrix space of
Dipper-Donkin is  an associative algebra  $D_q(n)$ generated by
$n^2$ generators $Z_{ij}$ subject to the defining relations:
\begin{eqnarray}Z_{st}Z_{ij}&=& Z_{ij}Z_{st}+(q^2-1)Z_{it}Z_{sj},
\text{ if }s>i,t>j,\\\nonumber Z_{st}Z_{ij}&=&q^{2}Z_{ij}Z_{st},
\text {if }s>i,t\le j,\\\nonumber
Z_{it}Z_{ij}&=&Z_{ij}Z_{it},\end{eqnarray}

Let $G$ be the semi-group ${\mathbb Z}_+^n$ with standard basis
$e_1, e_2,\cdots, e_n$. The algebra   $D_q(n)$ is $G\times
G$-graded with

$$deg Z_{ij}=(e_i, e_j).$$

Let $\phi$ and $\psi$ be semi-group homomorphisms defined by
\begin{eqnarray}\phi: G\times G &\longrightarrow & K^*, \\\nonumber
(e_i, e_j)&\mapsto & ¡¡p_{ij},  \text{ for all }i,j.\end{eqnarray}
and

\begin{eqnarray}\psi: G\times G &\longrightarrow & K^*, \\\nonumber
(e_i, e_j)&\mapsto & q_{ij}, \text{ for all }i,j.\end{eqnarray}

From the defining relations, one can see clearly that  the algebra
${\s O}_{q, P, Q}(M(n))$ is a bi-character deformation of the
algebra $D_q(n)$. \qed

\begin{Rem}The algebra ${\s O}_{q,P,Q}(M(n))$ is an iterated Ore
extension  and hence a noetherian domain. If we put
$p_{ij}=q_{ji}$ we get the algebras constructed by
Artin-Schelter-Tate. The algebras considered in \cite{jz2} are
also special cases of the present algebras.
\end{Rem}

Our construction of the basis depends heavily on a bar action. The
following can be obtained easily from the defining relations of
the algebra ${\s O}_{q,P,Q}(M(n))$.

\begin{Lem}The assignment
\begin{eqnarray}^{-}:{\s O}_{q,P,Q}(M(n))\longrightarrow
{\s O}_{q,P,Q}(M(n)),\\\nonumber Z_{ij}\mapsto Z_{ij},q\mapsto
q^{-1},p_{ij}\mapsto p_{ji},q_{ij}\mapsto q_{ji}\end{eqnarray}
extends to an
 algebra anti-automorphism over ${\mathbb Q}$.\end{Lem}

Two elements $x,y$ in the algebra ${\s O}_{q,P,Q}(M(n))$ are
called equivalent if there is a monomial $m$ of
$q,q^{-1},p_{ij},q_{ij}$ such that $x=my$. In this case, we write
$x\sim y$.

For any matrix $A=(a_{ij})\in M_n({\mathbb Z}_+)$,  we define the
monomial  $Z^A=\Pi Z_{ij}^{a_{ij}}$, where the factors are
arranged according to the lexicographic ordering.

Using Bergman's diamond lemma, we see that the algebra ${\s
O}_{q,P,Q}(M(n))$ has the nice basis

$$\{Z^A|A\in M_n({\mathbb Z}_+)\}.$$

However, the above basis is only almost the right choice for our
construction; another normalization will be needed:

Let
$$D(A)=q^{\sum_{s>i,j\le t}a_{st}a_{ij}} \Pi_{s>i}
p_{si}^{a_{st}a_{ij}}\Pi_{t>j,s\le
i}q_{tj}^{a_{st}a_{ij}}\Pi_{s>i,j>t}q_{jt}^{-a_{st}a_{ij}},$$
we
define the normalized monomial $Z(A)=D(A)Z^A$. For matrices $A$
and $B$ in $M_n({\mathbb Z}_+)$, we define  $B\le A$ if $B$
can be obtained from $A$ by a sequence of $2\times 2$ sub-matrix
moves  of the form

\begin{equation}\label{m-form}\begin{pmatrix}a_{ij}&a_{it}\\a_{sj}&a_{st}\end{pmatrix}\longrightarrow
\begin{pmatrix}a_{ij}-1&a_{it}+1\\a_{sj}+1&a_{st}-1,\end{pmatrix}\end{equation}

where $a_{ij},a_{st}\ge 1$. Denote by $c_i$ the sum of the
elements in the $i$th column and $r_i$ the sum of the elements in
the $i$th row. Notice that if $B\le A$, then $A, B$ have the same
row sums and column sums.

 From the defining relations of
the algebra we have
$$\overline{Z(A)}=Z(A)+\sum_Bc_{AB}Z(B)$$
where $c_{AB}\in {\mathbb Z}[q,q^{-1},p_{ij},q_{ij}]$ and
$c_{AB}\ne 0$ only if $B\le A$.

Let $\Gamma$ be the subgroup of $K^*$ generated by $q, p_{ij},
q_{ij}$ for all $i,j=1,2,\cdots,n$.  Define an ordering on the
monomials of parameters by  $q<p_{ij}<q_{st}<0$ for $i<j; s<t$, and
extend to a lexicographic ordering on $\Gamma$ which is
compatible with the group structure on $\Gamma$. Denote by
$\Gamma_+$ the set of strictly positive elements and $\Gamma_-$
the set of strictly negative elements. Clearly,
$\Gamma_-=(\Gamma_+)^{-1}$.

 Using
Lusztig's method in section 2, we get

\begin{Thm}
For each $A\in M_n({\mathbb Z}_+)$, there is a unique element
$b(A)$ characterized by the following properties:
\begin{enumerate}
\item $b(A)=Z(A)+\sum_{B}h_{AB}Z(B)$, where $h_{AB}\in {\mathbb
Z}[\Gamma_+] $ and $h_{AB}\ne 0$ only if $B\le A$.

\item $\overline{b(A)}=b(A)$.
\end{enumerate}
The set $B^*=\{b(A)|A\in M_n({\mathbb Z}_+)\}$ is a ${\mathbb
Q}(q, p_{ij},q_{ij}|i,j=1,2\cdots,n)$ basis of the algebra ${\s
O}_{q,P,Q}(M(n))$.\end{Thm}

\begin{Rem}Later on, after proving the exponential nature we
  shall see that the coefficients $h_{AB}$ are in fact
  polynomials of $q^{-1}$. This fact justifies the word
  canonical.\end{Rem}

\section{The case of $2\times 2$}

\medskip

The bar action on the monomials $Z(A)$ is simply a reordering of the
generators. The purpose of  introducing the normalized monomial is
that it allows us to ignore the quasi-polynomial moves.  Clearly, only the
first relation in the defining relations produces new terms in the
process of reordering the generators. However, this relation only
involves four generators. In other words, only a $2\times 2$
sub-matrix of the matrix $A$ is involved. Hence, the bar action
can somehow be computed locally which means we should first
consider the $2\times 2$ case. The coordinate algebra of $2\times
2$ quantum matrix is an algebra with four generators and
relations:

\begin{eqnarray}Z_{22}Z_{11}&=& p_{21}^2q_{21}^2Z_{11}Z_{22}+(q^2-1)p_{21}^2Z_{12}Z_{21},\label{4.1}
\\\nonumber
Z_{22}Z_{21}&=& q_{21}^2Z_{21}Z_{22},\\\nonumber Z_{12}Z_{11}&=&
q_{21}^2Z_{11}Z_{12}\\\nonumber Z_{22}Z_{12}&=&
q^2p_{21}^2Z_{12}Z_{22},\\\nonumber Z_{21}Z_{12}&=&
q^2p_{21}^2q_{21}^{-2}Z_{12}Z_{21},\\\nonumber Z_{21}Z_{11}&=&
q^2p_{21}^2Z_{11}Z_{21}.
\end{eqnarray}

For $A=(a_{ij})_{2\times 2}\in M_2({\mathbb Z}_+)$, set
$$Z^A=Z_{11}^{a_{11}}Z_{12}^{a_{12}}Z_{21}^{a_{21}}Z_{22}^{a_{22}}.$$
The set $\{Z^A|A=(a_{ij})_{2\times 2}\in M_2({\mathbb Z}_+)$ is a
basis of the algebra. Denote by $\tr(A)=a_{11}+a_{22}$ the trace of
$A$ and $\str A=a_{11}-a_{22}$ the super trace of $A$.

Define $$D(A)=q^{(a_{11}a_{21}+a_{12}a_{21}+a_{12}a_{22})}
p_{21}^{(a_{11}a_{21}+a_{11}a_{22}+a_{12}a_{21}+a_{12}a_{22})}
q_{21}^{(a_{11}{\color{red}a_{12}}+a_{11}a_{22}-a_{12}a_{21}+{\color{red}a_{21}}a_{22})}$$
and $$Z(A)=D(A)Z^A$$

The $2\times 2$ quantum determinant is
$${\det}_q=Z_{11}Z_{22}-q_{21}^{-2}Z_{12}Z_{21}$$
which satisfies

\begin{eqnarray}\label{eqns}
{\det}_q \cdot Z_{11}=p_{21}^2q_{21}^2 \cdot Z_{11}\cdot{\det}_q &,&
Z_{22}\cdot {\det}_q=p_{21}^2q_{21}^2\cdot{\det}_q \cdot Z_{22}\\ \nonumber
{\det}_q\cdot
Z_{12}=q^2p_{21}^2q_{21}^{-2}\cdot{\color{red}Z_{12}}\cdot{\det}_q &,&
Z_{21}\cdot{\det}_q=q^2p_{21}^2q_{21}^{-2}\cdot{\det}_q\cdot 
Z_{21}.\end{eqnarray}

The normalized determinant
\begin{eqnarray}\Delta&=& p_{21}q_{21}\cdot {\det}_q\\&=& p_{21}q_{21}\cdot Z_{11}Z_{22}-p_{21}q_{21}^{-1}\cdot Z_{12}Z_{21}=
Z\begin{pmatrix}1&0\\0&1\end{pmatrix}-q^{-1}
Z\begin{pmatrix}0&1\\1&0\end{pmatrix}\end{eqnarray} is more appropriate for our
computations.

Examples:
$$b\begin{pmatrix}2&0\\0&1\end{pmatrix}=Z\begin{pmatrix}2&0\\0&1\end{pmatrix}-q^{-2}
Z\begin{pmatrix}1&1\\1&0\end{pmatrix},$$
$$b\begin{pmatrix}1&0\\1&1\end{pmatrix}=Z\begin{pmatrix}1&0\\1&1\end{pmatrix}
-q^{-1}Z\begin{pmatrix}0&1\\2&0\end{pmatrix},$$
$$b\begin{pmatrix}1&1\\0&1\end{pmatrix}=Z\begin{pmatrix}1&1\\0&1\end{pmatrix}
-q^{-1}Z\begin{pmatrix}0&2\\1&0\end{pmatrix}.$$

\smallskip

In the sequel we let $E=E_{12}+E_{21}$. Notice that the quantities
$r_2-r_1=a_{21}+a_{22}-a_{11}-a_{12}$,
$c_2-c_1=a_{12}+a_{22}-a_{11}-a_{21}$, and $a_{21}-a_{12}$ are the
same for $A, A\pm I$, and $A\pm E$. Set $f_A=q^{a_{21}-a_{12}}p_{21}^{r_2-r_1} q_{21}^{c_2-c_1}$.

\smallskip

\begin{Lem}\label{l41}
$Z(A)\Delta= f_A\left(Z(A+I)-q^{-\tr(A)-1}Z(A+E)\right)$. 
\end{Lem}

\pof This follows by an elementary computation using (\ref{4.1}) and (\ref{eqns}).\qed

\bigskip
 The next result is proved by similar arguments:

\begin{Lem}\label{l42}
$\overline{\Delta}=\Delta$. Moreover, 
\begin{displaymath}
  Z^A\Delta=f_A^2\Delta Z^A.
\end{displaymath}
\end{Lem}

\bigskip

Some of the following is well-known (\cite{fad1}, \cite{fad2}, \cite{koorn}), but the connection to
canonical bases seems to be new. 

\medskip

\begin{Def}We set 
\begin{displaymath}
  (n)_q=1+q^2+\dots +q^{2n-2}  \textrm{ if }n\geq2; (1)_q=1,\textrm{
    and }(0)_q=1.  
\end{displaymath}
\end{Def}

\begin{Def}We set
\begin{displaymath}
  \exp_q(X)=\sum_{n=0}^\infty \frac{X^n}{(n)_q!}.
\end{displaymath}
\end{Def}
\begin{Lem}If $XY=q^2YX$, then\label{l1}
\begin{displaymath}
  (X+Y)^n=\sum_{i=0}^n\binom{n}{i}_qX^iY^{n-i},
\end{displaymath}
where the coefficients satisfy
\begin{displaymath}
  \binom{n+1}{i}_q=\binom{n}{i}_q +q^{2n-2i+2}\binom{n}{i-1}_q ,
\end{displaymath}
and are given by
\begin{displaymath}
  \binom{n}{i}_q=\frac{(n)_q!}{(i)_q(n-i)_q}.
\end{displaymath}
\end{Lem}
\medskip

\begin{Lem}\label{l2}These quantized binomial coefficients satisfy the
following identity as follows easily by induction:
\begin{displaymath}
  \sum_{m=0}^s(-1)^mq^{m(m-1)}\binom{s}{m}_q=0. \end{displaymath}
\end{Lem}

\begin{Prop}If $YX=q^2XY$, then
\begin{displaymath}
  \exp_q(X+Y)=\exp_q(X)\exp_q(Y).
\end{displaymath}
This follows easily from  Lemma~\ref{l1} above.
\end{Prop}
\medskip

\begin{Prop}
\begin{displaymath}
  \exp_q(X)\exp_{q^{-1}}(-X)=1.
\end{displaymath}
\end{Prop}

\pof This follows easily from Lemma~\ref{l2} above. \qed

\vskip1cm

If we introduce change-of-basis matrices between the canonical basis
and the PBW basis, 

\begin{eqnarray*}
  Z(A)&=&\sum_CT_{CA}b(C)\\
\overline{Z(A)}&=&\sum_Dh_{DA}Z(D)\\
&=&\sum_C\overline{T_{CA}}b(C),\quad\textrm{ thus,}\\
TH&=&\overline{T}, \textrm{ i.e.}\\
H&=&T^{-1}\overline{T}.
\end{eqnarray*}

\medskip

Define operators $t$ and $\overline{t}$ as follows: For  a
matrix $A=\begin{pmatrix}a_{11}&a_{12}\\a_{21}&a_{22}\end{pmatrix}\in
M_2({\mathbb Z}_+)$ with $a_{11}a_{22}>0$, define
$$A^\prime=\begin{pmatrix}a_{11}-1&a_{12}+1\\a_{21}+1&a_{22}-1\end{pmatrix}$$
The linear operator $t$ is given by
$$t(Z(A))=(\sum_{s={\mid} \str
A{\mid}+1}^{\tr(A)-1} q^{-s})Z(A^\prime)$$ if $a_{11}a_{22}>0$, and
zero otherwise. The linear operator $\overline{t}$ is given by
$$\overline{t}(Z(A))=(\sum_{s={\vert} \str
A{\vert}+1}^{\tr(A)-1} q^s) Z(A^\prime)$$ if $a_{11}a_{22}>0$, and
zero otherwise.

\medskip

The first result is straightforward and we omit the proof:

\begin{Lem}$\overline{t}t={q^{-2}}t\overline{t}$.
\end{Lem}

Let us define $\tau^{(i)}_{C}$  and
$\mu^{(i)}_{C}$  by 
\begin{eqnarray*}
t^i(Z(C))&=& \tau^{(i)}_{C}\cdot Z(C-iI+iE)\textrm{ and}\\
(-t+\overline{t})^i(Z(C))&=& \mu^{(i)}_{C}\cdot
Z(C-iI+iE),\end{eqnarray*} respectively. The first fundamental observation is

\begin{Lem}[Key]\label{key}
  \begin{eqnarray}\label{key1}
    \frac{\tau^{(r)}_{A+I}}{(r_{q^{-1}})!}&=&
\left(\frac{\tau^{(r)}_{A}}{(r_{q^{-1}})!}+ q^{-\tr(A-(r-1)I)-1} \frac{\tau^{(r-1)}_{A}}{((r-1)_{q^{-1}})!}\right),\\
\frac{\mu^{(s)}_{A+I}}{(s_{q^{-1}})!}&=&\label{4.6}\\\nonumber
\frac{\mu^{(s)}_{A}}{(s_{q^{-1}})!}&-& 
  q^{-\tr(A-(s-1)I)-1}
  \frac{\mu^{(s-1)}_{A}}{((s-1)_{q^{-1}})!}+
  q^{\tr(A)+1}  \frac{\mu^{(s-1)}_{A}} 
{((s-1)_{q^{-1}})!} \ .
  \end{eqnarray} 
\end{Lem}

\pof Let us set $A_{11}=a$ and $a_{22}=b$. We assume that $a\geq
b$. Then,
\begin{displaymath}
  tZ(A)=q^{-a+b-1}(b)_{q^{-1}}Z(A-I+E).
\end{displaymath}
Equation (\ref{key1}) follows from this and the simple identity
\begin{displaymath}
  (b+1)_{q^{-1}}=(b-r+1)_{q^{-1}}+(r)_{q^{-1}}q^{-2b+2r-2}.
\end{displaymath}
The second identity follows by a simple computation from the
observation that
\begin{displaymath}
 (\overline{t}-
 t)Z(A)=\frac{(q^{a}-q^{-a})(q^{b}-q^{-b})}{q-q^{-1}}Z(A-I+E). 
\end{displaymath}
\qed

\medskip

\begin{Prop}In the case at hand,
$$H=\exp_{{q^{-1}}}(-t+\overline{t}).$$
Moreover, $T=\exp_{q}(t)$, $\overline{T}=\exp_{{q^{-1}}}(\overline{t})$,
and $T^{-1}=\exp_{{q^{-1}}}(-t)$.
\end{Prop}

\pof The claim is that
\begin{equation}
   \overline{Z(A+I)}=\sum_s
   \frac{\mu^{(s)}_{A+I}}{(s_{q^{-1}})!}Z(A+I-sI+sE). 
\end{equation}

It follows easily from Lemma~\ref{l41} and Lemma~\ref{l42} that 
\begin{equation}
  \overline{Z(A+I)}=f_A^{-1}\overline{Z(A)}\Delta +q^{\tr(A)
    +1}\overline{Z(A+E)}. 
\end{equation}
Applying (\ref{4.6}), and Lemma~\ref{l41} once again, the claim follows
by induction on $\tr(A)$. \qed

\medskip

We have:

\begin{Lem}For any matrix $D\in M_n({\mathbb Z}_+)$,
\begin{displaymath}
  b(D)=\sum_C(T^{-1})_{C,D}Z(C).
\end{displaymath}
\end{Lem}

\pof It is clear that $\sum_C(T^{-1})_{C,D}Z(C)$ is bar invariant with 
the right leading term $Z(D)$. \qed

\medskip

\begin{Prop}
\begin{equation}
b(A)\cdot(f_A^{-1}\Delta)=b(A+I).
\end{equation}
\end{Prop}

\pof The claim is that
\begin{displaymath}
  \sum_r
  \frac{(-1)^r\tau^{(r)}_{A+I}}{(r_{q^{-1}})!}Z(A+I-rI+rE)=\sum_s\frac{(-1)^s\tau^{(s)}_{A}}{(s_{q^{-1}})!}Z(A-sI+sE) \cdot(f_A^{-1}\Delta).
\end{displaymath}
Using Lemma~\ref{l41}, this follows from (\ref{key1}) by induction
or $\tr(A)$. \qed

\bigskip

\section{The exponential nature}

\medskip

In this section, we compute the matrix of the bar action with
respect to the basis ${\mathbb B}_0$ consisting of normalized
monomials for general $n$. For $1\le i,j,s,t\le n$ with $i<s,j<t$
we define a linear operator $T_{ij}^{st}$ as follows:

\begin{eqnarray}T_{ij}^{st}(Z(A))=(\sum_{s=|a_{ij}-a_{st}|+1}^{a_{ij}+a_{st}-1} q ^{-s})
Z(A^\prime),\text{
if }a_{ij}a_{st}\ge 1,\\\nonumber T_{ij}^{st}(Z(A))=0,
\text{ if
}a_{ij}a_{st}=0.\end{eqnarray}

Similarly, an operator $\overline{T_{ij}^{st}}$ is defined as:

\begin{eqnarray}\overline{T_{ij}^{st}}(Z(A))=(\sum_{s=|a_{ij}-a_{st}|+1}^
{a_{ij}+a_{st}-1} q ^s)Z(A^\prime) \text{ if }a_{ij}a_{st}\ge
1,\\\nonumber \overline{T_{ij}^{st}}(Z(A))=0\text{ if
}a_{ij}a_{st}=0.\end{eqnarray}

The bar action  maps ${\mathbb B}_0$ into the ``ultimately
opposite'' basis ${\mathbb B}_u$ given by the elements
$\overline{Z(A)}$.

 We  know that for $n=2$, the matrix of the bar action is given by

${\mathbb H}_2=\exp_q(-T+\overline{T})$ as described in the last
section.

To facilitate the study of the bar action, we introduce a series
of intermediate PBW bases ${\mathbb
  B}_{k+1}=S_{i_k,j_k}^{s_k,t_k}{\mathbb B}_k$, where
$S_{i_k,j_k}^{s_k,t_k}$ is a linear map which is applied to each
of the vectors in the given basis. It is the map which sends any
{\bf normalized} monomial of the form $\dots Z_{i_k
j_k}^{a_{i_kj_k}} \dots   Z_{i_k,t_k}^{a_{i_kt_k}} \dots Z_{s_k
j_k}^{a_{s_kj_k}} \dots Z_{s_k t_k}^{a_{s_kt_k}}$ into the
corresponding normalized expression $\dots Z_{s_k
t_k}^{a_{s_kt_k}} \dots   Z_{s_k,j_k}^{a_{s_kj_k}} \dots Z_{i_k
t_k}^{a_{i_kt_k}} \dots Z_{i_k j_k}^{a_{i_kj_k}}$.  If we assume,
and this condition will always be satisfied in the applications,
that the interior ellipses represent terms which quasi-commute
with the elements $Z_{i_k
  j_k}, Z_{i_k,t_k},  Z_{s_k j_k} $, and $Z_{s_k t_k}$, then the
matrix  ${\mathbb S}_{i_k,j_k}^{s_k,t_k}$ of the map
$S_{i_k,j_k}^{s_k,t_k}$ with respect to the basis ${\mathbb B}_k$
is precisely ${\mathbb H}_2$ tensored appropriately with an
identity operator representing all the variables which stay fixed.
${\mathbb
  S}_{i_k,j_k}^{s_k,t_k}$ is also the change of basis matrix from
${\mathbb B}_{k+1}$ to ${\mathbb B}_k$.

To give the matrix ${\mathbb H}$ of the bar action in the basis
${\mathbb B}_0$ is the same as giving the change of basis matrix
from ${\mathbb B}_u$ to ${\mathbb B}_0$.

\begin{Thm}\label{thm5.1}For any matrix $A\in M_n({\mathbb Z}_+)$,
\begin{equation}{\mathbb H} = \Pi_{(i,j),(s,t),i<s,j<t}
{\mathbb S}_{ij}^{st}\end{equation} where the factors are arranged
according to the double lexicographic ordering which means that we
first use lexicographic ordering on the indices $(i,j)$ and then
lexicographic ordering of indices $(s,t)$. The matrices ${\mathbb
  S}_{ij}^{st}$ are viewed as above.
\end{Thm}

\pof One only needs to verify that at each step, the ellipses
mentioned in the discussion above indeed do represent
quasi-commuting element. This is elementary. \qed

\begin{Rem}There are at least three other such decomposition, namely
  where one uses the opposite ordering in one or both places.
\end{Rem}

Theorem~\ref{thm5.1} shows that the bar action on the normalized
monomials only depends on $q$, so the canonical basis only depends on
the parameter $q$. This also means that for the multi-parameter case,
the expression of the dual canonical basis is exactly the same as that
expression in the one-parameter case. In \cite{zhang}, it is proved
that the dual canonical basis is invariant under the multiplications
of certain covariant quantum minors which is also true by the above
theorem.

Letting  $p_{ij}=q_{ij}$, we get the algebra constructed in
\cite{ast} which is a bi-algebra with the usual coproduct.
Inverting the quantum determinant we get the quantum function
algebra which is dual to the quantum enveloping algebra $U_{q,
P}(gl_n)$ (see \cite{ast} for more detail), as well as a basis
$\bar{B}$ of this quantum function algebra. The dual basis of
$\bar{B}$ is a basis of $U_{q, P}(gl_n)$ and this basis also only
depends on the parameter $q$.

The upper triangular case enable us to construct a basis of
$U_q(n^+)$ which should be the canonical basis constructed by
Lusztig. To this end we need to show that the basis consisting of
the images of $Z(A)$ is dual to the PBW basis consisting of
divided powers.

\section{An inductive program}

We only need to consider the official one. For any matrices $A,
B\in M_n({\mathbb Z}_+)$, there exist $d_{AB}, d_{BA}\in{\mathbb
Z}$ such that
$$Z(A)Z(B)=q^{d_{AB}}Z(A+B)+\text{ lower order terms}$$
and
$$Z(B)Z(A)=q^{d_{BA}}Z(A+B)+\text{ lower order terms}.$$

A direct computation shows that $d_{AB}=-d_{BA}$.

Assume that both $b(A)$ and $b(B)$ are known, then
$\frac{qb(A)b(B)-q^{-1}b(B)b(A)}{q-q^{-1}}$ is invariant under the
bar action and so

$\frac{q^{1-d_{A,B}}b(A)b(B)-q^{d_{A,B}-1}b(B)b(A)}{q-q^{-1}}=
b(A+B)+\sum_{D<A+B}c_{A,B}^Db(D) $ with the coefficients
$c_{A,B}^D\in {\mathbb Z}[q+q^{-1}]$. Hence the element $b(A+B)$
can be determined uniquely by the above equation.

After showing the exponential nature, we see that the expression
of the dual canonical basis elements do not depend on the choice
of the parameters. Hence, the dual canonical basis is stable 
(up to the equivalence relation) under multiplication by 
covariant minors.

Now, let us compute the basis for $2\times 3$ matrices. 
After the removal of the covariant minors,  the only
case we need to compute directly is for
the sub-matrices

$$\begin{pmatrix}a&0&0\\0&b&c\end{pmatrix}$$

\begin{enumerate}
\item If $a\le b$, the basis element is
$$q^{-bc}(Z_{11}Z_{22}-q^2Z_{12}Z_{21})^a Z_{22}^{b-a}Z_{23}^c$$.

\item If $a>b$ and $a\ge b+c$, then the basis element is
$$q^{-bc}Z_{11}^{a-b-c}(Z_{11}Z_{22}-q^2Z_{12}Z_{21})^b(Z_{11}Z_{23}-q^2Z_{13}Z_{21})^c$$.

\item if $a>b$ but $a<b+c$, then the basis element is
$$q^{-bc}(Z_{11}Z_{22}-q^2Z_{12}Z_{21})^b(Z_{11}Z_{23}-
q^2Z_{13}Z_{21})^{a-b}Z_{23}^{c-a+b}$$.
\end{enumerate}

\medskip

\section{ A conjecture}

Given a matrix $A\in M_n({\mathbb Z}_+)$ 
we can draw a graph $H(A)$. The nodes of the graph are the
matrices $B$ which can be obtained from $A$ by $2\times 2$ matrix
transformations of the form (\ref{m-form}).
 We place $B_1$ on a level above $B_2$ if
$B_2$ can be obtained from $B_1$. We draw a line between two
nodes if we get the lower node from the upper one by one
transformation; we  attach to the line the upper indices $(i,j)$
and lower indices $(s,t)$ with the obvious meaning. Notice that, like
at the top,  there
is just one graph at the bottom of the graph,   say $T_A$; the tail of $A$.

A path in $H(A)$ is called principal path from $A$ to $B$ if it is
a longest path from $A$ to $B$ and is maximal (according to the
lexicographic order) among all of the longest paths in $H(A)$
from $A$ to $B$.

  We use $p_{ij}^{st}(B)$ to denote the
number of lines  in the principal path from $A$ and $B$  with
indices $(i,j),(s,t)$. Denote by $l_{AB}$ the number of lines (
the length) from $A$ to $B$ in the principal path.

\begin{Con}
$$b(A)=\sum_{B\le
A}(- q ^{-1})^{l_{AB}}\Pi_{(i,j),(s,t)i<s,j<t} q
^{-|a_{ij}-a_{st}|}
\begin{pmatrix}min \{ a_{ij},a_{st}\}\\ p_{ij}^{st}(B)\end{pmatrix}_{ q ^{-2}}Z(B).$$
\end{Con}

We only need to prove that the elements are bar invariant. A
direct computation shows that the conjecture hold for the cases
$2\times 2$ and $2\times 3$.

\section{The positivity}

\medskip

The negative part $U_q(n^-)$ is the subalgebra of quantum
enveloping algebra $U_q(A_{2n-1})$ generated by  $F_1, F_2,
\cdots, F_{2n-1}$  subject to the quantum Serre relations:

\begin{eqnarray}F_iF_j=F_jF_i, \text{ if }|i-j|>1,\\\nonumber
F_i^2F_j-(q^2+q^{-2})F_iF_jF_i+F_jF_i^2=0, \text{ if
}|i-j|=1.\end{eqnarray}

 For a homogeneous element $x$, denote by $wt(x)$ the weight of
$x$.

   Let $\Pi=\{\alpha_1, \alpha_2,
\cdots, \alpha_{2n-1}\}$ be the set of simple roots of the Lie
algebra of type $A_{2n-1}$. For a positive root
$\alpha_{ij}=\alpha_i+\alpha_{i+1}+\cdots+\alpha_j$, the
corresponding root vector is defined as
$$F_{\alpha_i}=F_i, F_{ij}=F_{\alpha_{ij}}:=[F_j,[F_{j-1},[\cdots,[F_{i+1}, F_i]_q\cdots
]_q]_q]_q. $$
The quantum commutator is given by
$$[x,y]_{q}=xy-q^{-2(\alpha, \beta)}yx,$$
for homogeneous elements $x$ and $y$ with weights $\alpha$ and
$\beta$ respectively.

We introduce an ordering on the set of positive roots  $\Delta_+$ by
$$\alpha_{ij}\le\alpha_{kl}, \text{ if }j<l\text{ or }j=l\text{
and } i<k.$$ Notice that this ordering coincides with the ordering given
by the reduced expression of the longest element $w_0=r_1r_2\cdots
r_{2n-1}r_2\cdots r_{2n-1}\cdots r_{2n-2}r_{2n-1}r_1$ in the Weyl
group, where $r_i$ is the simple reflection determined by the
simple root $\alpha_i$. The PBW basis of $U_q(n^-)$ is indexed by
the set ${\mathbb Z}_+^{\Delta_+}$.

For an $\bold{m}=(m_{ij})\in {\mathbb Z}_+^{\Delta_+}$, denote by
$\Deg{\bold m}=\sum m_{ij}(j-i+1)$. The PBW basis element indexed by
$\bold{m}$ is
$$F({\bold m}):=\Pi\frac{F_{ij}^{m_{ij}}}{[m_{ij}]_{q^2}!},$$
where the factors are arranged according to the above ordering on the
set $\Delta_+$. Denote by $|{\bold m}|$ the weight of $F({\bold
m})$.

The tensor product $U_q(n^-)\otimes U_q(n^-)$ can be regarded as a
${\mathbb Q}(q)$-algebra with multiplication
$$(x_1\otimes x_2)(x_1^\prime\otimes
x_2^\prime)=q^{2(wt(x_2),wt(x_1^\prime))}x_1x_1^\prime\otimes
x_2x_2^\prime,$$ for homogeneous elements $x_1,x_2,x_1^\prime,
x_2^\prime\in U_q(n^-)$.

In \cite{lu3}, it was proved that
\begin{Lem}The following assignment
\begin{eqnarray}r:U_q(n^-)\longrightarrow U_q(n^-)\otimes
U_q(n^-)\\\nonumber F_i\mapsto F_i\otimes 1+1\otimes F_i,\text{
for all }i\end{eqnarray} extends to an algebra
homomorphism.\end{Lem}

\begin{Rem}The algebra homomorphism $r$ is
co-associative.\end{Rem}

There is a scalar product on $U_q(n^-)$ (see \cite{lu3})
satisfying
$$ (F_i,
F_j)=\delta_{ij}, (x, y_1y_2)=(\Delta(x), y_1\otimes y_2),
(x_1x_2, y)=(x_1\otimes x_2, \Delta(y)),$$ where the scalar
product on $U_q(n^-)\otimes U_q(n^-)$ is given by
$$(x_1\otimes x_2,y_1\otimes
y_2)=(x_1, y_1)(x_2, y_2).$$ On the PBW basis, the scalar product
is given by
$$(F(\bold{m}), F(\bold{n}))=\frac{(1-q^4)^{\Deg{\bold
m}}}{\Pi_{i\le j}\phi_{m_{ij}}(q^4)}\delta_{{\bold m}, {\bold
n}},$$ where $\phi_k(z)=(1-z)(1-z^2)\cdots (1-z^k)$.

Let
$${\mathcal L}:=\oplus_{m\in {\mathbb Z}_+^{\Delta_+}}
 {\mathbb Z}[q]F({\bold m}).$$

 Denote by $-$ the ring automorphism:

\begin{eqnarray} -: U_q(n^-)\longrightarrow U_q(n^-),\\\nonumber
\overline{F_i}=F_i, \overline{q}=q^{-1}\text{ for all
}i.\end{eqnarray}

The canonical basis $B=\{G({\bold m})\}$ (the lower global crystal
basis in Kashiwara's terminology) of $U_q(n^-)$  is a  ${\mathbb
Z}[q]$ basis of ${\mathcal L}$ such that

$$\overline{G({\bold m})}=G({\bold m}), G({\bold m})=F({\bold m})\quad
\mod q{\mathcal L}.$$

Lusztig proved that the canonical basis enjoys some remarkable
properties:

\begin{Thm}\label{8.3} (Positivity, \cite{lu3})The following hold.
\begin{enumerate}
\item For any $b, b^\prime\in B$, we have
$$bb^\prime=\sum_{b''\in B, n\in{\mathbb Z}}
c_{b,, b^\prime, b'',n}q^nb'',$$ where $c_{b,, b^\prime, b'',n}\in
{\mathbb Z}_+$ are zero except for finitely many $b'', n$.

\item For any $b\in B$, we have
$$r(b)=\sum_{b',b''\in B, n\in{\mathbb Z}}d_{b, b', b'',n}q^n
b'\otimes b'',$$ where $d_{b, b', b'',n}\in{\mathbb Z}_+$ are zero
except for finitely many $b', b'', n$.

\item For any $b, b^\prime\in B$ we have

$$(b, b^\prime)=\sum_{n\in{\mathbb Z}_+}f_{b, b^\prime, n}q^n$$
where $f_{b, b^\prime, n}\in {\mathbb Z}_+$.
\end{enumerate}
\end{Thm}

 The canonical basis $B$ is almost orthogonal with respect to the above scalar product.
 By almost orthogonal one means that

$$(G({\bold m}), G({\bold n}))=\delta_{{\bold m}{\bold n}}\mod q{\mathcal A}, $$
where ${\mathcal A}$ is the subring of ${\mathbb Q}(q)$ consisting of
the rational functions regular at $q=0$.

Denote by $\{F^*({\bold m})\}$ and $\{G^*({\bold m})\}$ the dual
bases of $\{F({\bold m})\}$ and $\{G({\bold m})\}$ respectively
with respect to the above scalar product. Since $\{F({\bold m})\}$
is orthogonal, the basis $\{F^*({\bold m})\}$ is simply a
rescaling  of $\{F({\bold m})\}$, namely

$$F^*({\bold m})=\Pi q^{\binom{m_{ij}}{2}}
{F^*_{ij}}^{m_{ij}},$$ where 
$$F^*_{ij}= (1-q^4)^{i-j}F_{ij}.$$

The dual canonical basis $\{G^*({\bold m})\}$ can be characterized
by conditions similar to those defining the canonical basis. Let

$${\mathcal L}^*:=\oplus_{\bold m} {\mathbb Z}[q]F^*({\bold m}).$$

  Let $\Phi$ be the anti-automorphism of
$U_q(n^-)$ such that $\Phi(F_i)=F_i$, and $\Phi(q)=q^{-1}$.
 In \cite{lnt}, it was proved that

\begin{Prop}Let ${\bold m}\in {\mathbb Z}_+^{\Delta_+}$
and write $|{\bold m}|^2:=(wt({\bold m}), wt({\bold m}))$. Then
$G^*({\bold m})$ is the unique homogeneous element of degree
$wt({\bold m})$ of $U_q(n^-)$ satisfying

$$\Phi(G^*({\bold m}))=q^{2\Deg{\bold m}-|{\bold m}|^2}G^*({\bold m}), G^*({\bold m})=F^*({\bold m})\quad
\mod q{\mathcal L}^*.$$\end{Prop}

The dual canonical basis of $U_q(n^-)$ can be constructed using
Lusztig's elementary method by modifying the above construction.
Let $F^*_N({\bold m})=q^{\frac{1}{2}|{\bold
    m}|^2-\Deg{\bold m}}F^*({\bold m})$ and let $G^*_N({\bold
  m})=q^{\frac{1}{2}|{\bold m}|^2-\Deg{\bold m}}G^*({\bold m})$. Let

$${\mathcal L}^*_N=\oplus_m{\mathbb Z}[q]F^*_N({\bold m}).$$

Then the above proposition can be rewritten as

\begin{Prop}Let ${\bold m}\in {\mathbb Z}_+^{\Delta_+}$. Then
$G^*_N({\bold m})$ is the unique homogeneous element of degree
$wt({\bold m})$ of $U_q(n^-)$ satisfying

$$\Phi(G^*_N({\bold m}))=G^*_N({\bold m}), G^*_N({\bold m})=F^*_N({\bold m})\quad
\mod q{\mathcal L}^*_N.$$\end{Prop}

The basis $\{G^*_N({\bold m})|{\bold m}\in {\mathbb
Z}_+^{\Delta_+}\}$ is called the normalized dual canonical basis
of $U_q(n^-)$.

\begin{Thm} The multiplication of  the basis $B^*$
of the algebra $O_q(M(n))$ satisfies a positivity property
analogous to 1 in Theorem~\ref{8.3}.\end{Thm}

\pof By duality, the positivity properties hold for the dual
canonical basis of the canonical basis of $U_q(A_{2n-1})^-$.
Hence, The multiplication of  the normalized dual canonical basis
$\{G^*_N({\bold m})|{\bold m}\in {\mathbb Z}_+^{\Delta_+}\}$
 has the positivity property.

Define inductively

\begin{eqnarray}F_{i, i}^*&=&F_i\;,\qquad\textrm{ and}\\
 F_{i, j}^*&=&\frac{qF_jF_{i,j-1}^*-q^{-1}F_{i,j-1}^*F_j}{q^2-q^{-2}}.
\end{eqnarray}
A simple induction on $j-i$ gives that the $F_{ij}^*$'s are fixed by
$\Phi$.  It is easy to see that the $F_{ij}^*$ satisfy the defining
equations of the quantum matrix indeed, that one may embed the algebra
$O_q(M(n))$ into $U_q(n^-)$ as a subalgebra by:
$$Z_{ij}\mapsto F_{i,j+n}^*\text{ for all }i,j=1,2,\cdots, n.$$
Under this embedding, the PBW basis $\{Z(A)|A\in M_n({\mathbb
Z}_+)\}$ is a subset of the basis $\{F^*_N({\bold m})\}$ By the
descriptions of the basis $B^*$ obtained from Lusztig's procedure,
one sees that $B^*$ is a subset of the normalized dual canonical
basis of $U_q(n^-)$. \qed

\begin{Def}An element $b\in B^*$ is called
decomposable if there exist $m\in {\mathbb Z}$ and $b_1, b_2\in
B^*$ such that $b=q^m b_1b_2$. The basis element is called
indecomposable otherwise.\end{Def}

By definition, any basis element $b\in B^*$ can be written as

\begin{equation}\label{pos}b=q^m b_1b_2\cdots b_s\end{equation}
for some $m\in {\mathbb Z}$ and
indecomposables $b_1, b_2,\cdots, b_s\in B^*$.  Notice
  that it follows from the positivity that any $b$ written as in
  (\ref{pos}) as a
  product of more than two basis elements is, indeed, decomposable.

\begin{Con} The decomposition is unique up to a permutation.
\end{Con}

\medskip

\begin{Prop} Let $b\in B^*$. Assume that $b=q^m b_1b_2\cdots b_s$.
Then for any $i,j\in{\overline[1, n]}$ there exists an integer
$h_{ij}$ such that $b_ib_j=q^{h_{ij}}b_jb_i$. Furthermore, for any
$\{i_1<i_2<\cdots<i_r\}\subset{\overline[1, n]}$ there exists an
element $e\in B^*$ such that
$$e=q^sb_{i_1}b_{i_2}\cdots b_{i_r}.$$

\end{Prop}

\pof We use induction on $s$. The case of $s=1$ is trivial. Assume
that our hypothesis holds for $\leq s-1$. We first prove that $b_i$
$q$-commutes with $b_{i+1}$. By the positivity of the
multiplication $b_ib_{i+1}$ must be of the form $a_{ij}b'$ for
some $a_{ij}\in {\mathbb Z}_+[q, q^{-1}]$ and $a_{ij}$ must be a
power of $q$ since multiplication with basis elements can not
decrease the number of summands. Applying the bar action, we see
that
$$b_jb_i=\bar{a_{ij}}b'=\bar{a_{ij}}a_{ij}^{-1}b_ib_j.$$
Hence, all of
the factors $b_i$ and $b_j$ $q$-commute with each other. Because
multiplication satisfies positivity, the statement follows. \qed.

\smallskip

Hence, to understand the basis $B^*$ completely, one must determine
the indecomposables.

\medskip

\begin{Prop} All  quantum minors are
indecomposable.\end{Prop}

\pof Let $D$ be a quantum minor. Specifically, $D$ is the quantum
determinant of a subalgebra $A$ of $O_q(M(n))$ isomorphic to
$O_q(M(m))$ for some $m\le n$. If $D$ is decomposable and is written
as a product of two basis elements $b_1, b_2$, then $b_1, b_2$ are
members of the dual canonical basis of $A$. Hence, we may only
treat the case that $D$ is the quantum determinant.  With respect to
the lexicographic order, the leading term of $det_q$ is
$Z_{11}Z_{22}\cdots Z_{nn}$. Hence, the leading terms of $b_1$ and
$b_2$ produce the term $Z_{11}Z_{22}\cdots Z_{nn}$. But then the
leading term of $b_1$ will be of the form $Z_{i_1,i_1}\dots
Z_{i_r,i_r}; i_1<\cdots<i_r$, and the leading term of $b_2$ will be
what remains. But then $b_1$ is a minor centered around
the diagonal, and so is $b_2$. It is clear that the product of two
such ``disjoint'' minors cannot give the full quantum
determinant. 
\qed

\medskip

In \cite{z}, it was proved that the dual canonical basis $B^*$ of
the algebra $O_q(M(n))$ is invariant under multiplication
  by  the quantum determinant. Setting the quantum determinant to $1$, we
get a basis $K^*$ of the algebra $O_q(SL(n))$.  Clearly, we have

\begin{Thm} The multiplication of the basis $K^*$ of the algebra $O_q(SL(n))$
has  the positivity property.\end{Thm}

In \cite{z}, it was proved that the basis $K^*$ is dual to the
canonical basis of the modified quantum enveloping algebra
$\widetilde{U_q(A_{2n-1})}$ (one can refer to \cite{lu2} for more
details of the construction of the canonical basis of the modified
quantum enveloping algebra). By duality again, we have

\begin{Thm}The co-product of the canonical basis of
$\widetilde{U_q(A_{2n-1})}$ has the positivity property.\end{Thm}

This result was originally conjectured by Lusztig \cite{lu2}. We
believe that  this result holds  in all simply-laced cases. We will
deal with this  in a forthcoming paper.

\bibliographystyle{amsplain}

\end{document}